\documentclass{article}

\usepackage[english]{babel}

\usepackage{float}
\usepackage{geometry} 
\usepackage{amsmath}
\numberwithin{equation}{section}
\usepackage[ruled]{algorithm2e}
\usepackage{graphicx}
\usepackage[colorlinks=true, allcolors=blue]{hyperref}

\usepackage{xcolor}

\usepackage{authblk}

\usepackage{csquotes}

\usepackage{enumitem}

\usepackage{amsthm}
\usepackage{amssymb}
\usepackage{enumitem}

\newtheorem{theorem}{Theorem}[section]

\newtheorem{lemma}[theorem]{Lemma}

\newtheorem{example}{Example}[section]
\newtheorem{corollary}{Corollary}[section]

\newtheorem{remark}[theorem]{Remark}
\newtheorem{problem}{Problem}[section]

\title{A two-point  phase recovering with spherical wave reference}

\date{}

\author{R.G. Novikov, V.N. Sivkin}

\begin{document}
\maketitle

\begin{abstract}
    We consider a reference wave, a radiation solution, and the sum of these solutions (total solution) for the Helmholtz equation in an exterior region. We give two-point formulas for  approximate phase recovering of the radiation solution from the intensity of the total solution for the case of spherical reference wave. We show that these formulas can be used, in particular, for approximate phase recovering from holographic data on a single plane. By these formulas, we continue previous studies with plane wave reference.
\end{abstract}



\section{Introduction}


\noindent We consider the Helmholtz equation in an exterior region 
\begin{align}\label{eq:schrod}
    -\Delta \psi(x) = \kappa^2\psi(x),  \quad x \in  {\cal U}, \quad \kappa >0, \quad d\geq 2, 
\end{align}
where ${\cal U}$ is an exterior region in $\mathbb{R}^d,$ e.g.,  
\begin{align}\label{eq:1.2}
{\cal U} = \{x \in \mathbb{R}^d,\,\, |x|>\rho\}, \quad\rho\geq 0.     
\end{align}
For this equation we consider its solutions $\psi$ satisfying the Sommerfeld radiation condition
\begin{align}\label{eq:sommerf}
    &|x|^{\frac{d-1}{2}}\left(\frac{\partial}{\partial |x|} - i \kappa\right) \psi(x) \to 0 \quad \text{ as } |x| \to +\infty \quad \text{(uniformly in }\, x/|x|);
\end{align}
see, for example, \cite{Karp61}, \cite{W}.
Recall that
\begin{align}\label{eq:psif}
    \psi(x) = \frac{e^{i\kappa|x|}}{|x|^{\frac{d-1}{2}}}f\left(\frac{x}{|x|}\right)+ {\cal O}\left(\frac{1}{|x|^{\frac{d+1}{2}}}\right), \quad \text{ as } |x| \to +\infty,
\end{align}
where $f$ arising in \eqref{eq:psif} is known as the far-field pattern of $\psi.$

We also consider a spherical wave $\psi_0$ generated by a point source at $x_0 \in {\mathbb R^d}$, i.e., 
\begin{align}\label{eq:1.4}
\begin{aligned}
&(\Delta+\kappa^2) \psi_0(x) = A\delta(x-x_0), \quad x \in \mathbb{R}^d, \quad A \in \mathbb{C}\setminus \{0\},\\
&\psi_0 \text{ satisfies \eqref{eq:sommerf},}
\end{aligned}
\end{align}
where  $\delta$ is the Dirac function.

We consider the following problem.
\begin{problem}\label{prbl:1}
   Let $\psi_0$ be as in \eqref{eq:1.4}, and $\psi_1$ satisfy \eqref{eq:schrod}, \eqref{eq:sommerf}. Find $\psi_1$ on ${\cal U}'$ from  $|\psi_0+\psi_1|^2$ on ${\cal U}'',$ at fixed $\kappa,$ $A,$ $x_0,$  where   ${\cal U}',$ ${\cal U}''$  are some subsets of ${\cal U}.$ 
\end{problem}

This problem is motivated by holography and phaseless inverse scattering; see, e.g., \cite{G}, \cite{G49}, \cite{K2014}, \cite{N2015}.

Adopting holographic terminology, we say that $\psi_0$ is a reference beam, $\psi_1$ is an object beam, $\psi = \psi_0+\psi_1$ is the total field, $|\psi|^2$ on  ${\cal U}''$ is a hologram. 
Basically, in holography, ${\cal U}' = {\cal U}'' = X,$ where $X$ is a hyperplane in ${\cal U}.$

In general, Problem \ref{prbl:1} is not uniquely solvable even when ${\cal U}'' = {\cal U}.$ A simple example is as follows.

\begin{example}\label{ex:1}
 Let $\psi_0$ be as in \eqref{eq:1.4}, $x_0 \in \mathbb{R}^d \setminus{\cal U},$ where ${\cal U}$ is as in \eqref{eq:1.2}. Let   
\begin{align}
\begin{aligned}
&(\Delta + \kappa^2) \psi_1(x) = c\delta(x-x_0), \quad x \in \mathbb{R}^d, \quad c\in \mathbb{C}\setminus\{0\}, \\
&\psi_1 \text{ satisfies \eqref{eq:sommerf}, }
\end{aligned} 
\end{align}
where 
\begin{align}
  |A+c|^2 = |A|^2. \label{eq:Ac}  
\end{align}
Let ${\cal U}',$ ${\cal U}''$ be subsets of ${\cal U},$ where ${\cal U}'$ is such that
\begin{align}\label{eq:1.8}
   \text{ if $\psi$ satisfies \eqref{eq:schrod}--\eqref{eq:sommerf}, and $\psi \equiv 0$ on ${\cal U}',$ then $\psi \equiv 0,$ on }{\cal U}.
\end{align}

Then
\begin{align}
&|\psi_0+\psi_1|^2 \equiv |\psi_0|^2 \text{ on } {\cal U}'', \label{eq:6}\\
&\psi_1\not\equiv 0 \text{ on } {\cal U}'. \label{eq:7}   
\end{align}

\end{example}

In particular, in Example \ref{ex:1}, one can assume that ${\cal U}'' = {\cal U},$ ${\cal U}' = X,$ where $X$ is a hyperplane in ${\cal U}.$ 
Note also that condition \eqref{eq:Ac} is an equation of a circle with respect to $c.$


Note that property \eqref{eq:6} follows from the formula 
\begin{align}
\psi_0+\psi_1 = \frac{A+c}{A} \psi_0
\end{align}
and condition \eqref{eq:Ac}. In addition,  for \eqref{eq:7} we also use  assumption \eqref{eq:1.8}. 

On the other hand, for Problem \ref{prbl:1}, inspite of Example \ref{ex:1}, for large classes of  $x_0$ and $x,$ $y \in {\cal U}, $ there are explicit approximate formulas for finding $\psi_1(x)$ from the intensities $|\psi_0(x)+\psi_1(x)|^2$ and $|\psi_0(y)+\psi_1(y)|^2,$  and given parameters $A$ and $x_0$ in \eqref{eq:1.4}. Some  simple formulas of this type are given in Theorem \ref{thm:main} in Section \ref{sec:formulas}. These formulas are given in terms of the far-field pattern $f_1$ for $\psi_1,$ where $f_1$ is defined in  \eqref{eq:psif1}. In these formulas $|x|$ and  $|x_0|$ are sufficiently large. By these formulas, we continue previous studies with plane wave reference $\psi_0$; see \cite{N2015}, \cite{NSh}, \cite{NS26},   and references therein.
In addition, our formulas admit direct applications to Problem \ref{prbl:1} for the case when ${\cal U}' = {\cal U}'' = X,$ where $X$ is a hyperplane in ${\cal U},$ and the distance between $X$ and the origin $\{0\}\in \mathbb{R}^d$ is sufficiently large; see Corollary \ref{cor:1}, Lemma \ref{lem:2}, and Remark \ref{rem:3}   in Section \ref{sec:formulas}.

In larger framework, we contribute to holography and phaseless inverse scattering. These studies go back, in particular, to the pioneering works \cite{G}, \cite{G49}, \cite{W69}, \cite{W70}, and, for example,  more recent works \cite{D}, \cite{JL}, \cite{K2014}, \cite{M}, \cite{N2015}, \cite{Nu}. In connection with rather recent results in these directions, see, for example,  \cite{HNS},    \cite{ML}, \cite{NN}, \cite{N24}, \cite{NS26}, \cite{NS2026}, \cite{R}, and references therein. In fact, in the present work we continue studies of \cite{N2015} and very recent articles \cite{NS26}, \cite{NS2026}.  

The main results of this article are presented in details and proved in Sections \ref{sec:formulas}, \ref{sec:3}, \ref{sec:4}.

\section{Two-point formulas for phase recovering}\label{sec:formulas}

Let  $\psi_0$ be as in  \eqref{eq:1.4}. Then 
\begin{align}
&\psi_0(x) =  A G^+(x-x_0, \kappa), \\
&G^+(x, \kappa) =  -(2\pi)^{-d} \int_{\mathbb{R}^d} \frac{e^{i\xi x}d\xi}{\xi^2-\kappa^2-i0}, \quad x\in \mathbb{R}^d.
\end{align}
\noindent Note that, $G^+$ is the Green function for the Helmholtz operator $\Delta+\kappa^2$ with the Sommerfeld radiation condition \eqref{eq:sommerf}, and 
\begin{align}
&G^+(x, k) = -\frac{i}{4} H^1_0(|x||k|), \text{ for } d=2, \quad G^+(x, k) = -\frac{e^{i|k||x|}}{4\pi|x|}, \text{ for } d=3,
\end{align}
where $H^1_0$ is the Hankel function of the first type.

 Note  that 
\begin{align}
&\psi_0(x) = A' \frac{e^{i\kappa|x-x_0|}}{|x-x_0|^{(d-1)/2}}  +{\cal O}\left(\frac{A'}{|x-x_0|^{(d+1)/2}}\right), \text{ as } |x-x_0| \to +\infty, \label{eq:psi_green}\\
&A' = |A'|e^{i\varphi}:= A \frac{c(d, \kappa)}{(2\pi)^d}, \quad c(d, \kappa) = -\pi i (-2\pi i)^{(d-1)/2}\kappa^{(d-3)/2}. \label{eq:2.5}
\end{align}

 Let 
\begin{align}
\begin{aligned}\label{eq:b}
b(x, x_0) &:=\frac{|x|^{(d-1)/2}|x-x_0|^{(d-1)/2}}{|A'|}\left(|\psi_0+\psi_1|^2 - |\psi_0|^2\right), \quad x, x_0 \in \mathbb{R}^d, 
\end{aligned}   
\end{align}
where $\psi_0,$ $\psi_1$ are as in Problem \ref{prbl:1}, and $b$ depends also on $\psi_0,$ $\psi_1.$

In view of \eqref{eq:psif}, we have that
\begin{align}\label{eq:psif1}
    \psi_1(x) = \frac{e^{i\kappa|x|}}{|x|^{\frac{d-1}{2}}}f_1\left(\hat{x}\right)+ {\cal O}\left(\frac{1}{|x|^{\frac{d+1}{2}}}\right), \quad \hat{x} = \frac{x}{|x|}, \quad \text{ as } |x| \to +\infty,
\end{align}
where $f_1$ is the far-field pattern of $\psi_1.$

We consider $x,$ $x_0,$ $y$ such that 
\begin{align}
\begin{aligned}\label{eq:notations}
&x = r\hat{x}, \quad x_0 = \alpha r \hat{x}_0, \quad y = x+\zeta, \quad r>0, \\
&\hat{x},\, \hat{x}_0 \in \mathbb{S}^{d-1}, \quad \zeta \in \mathbb{R}^d, \quad \alpha>0,
\end{aligned}
\end{align}
where $\hat{x},$ $\hat{x}_0,$ $\zeta,$ and $\alpha$ are fixed, whereas $r$ is a parameter, which is sufficiently large, so that $x,$  $y\in {\cal U}.$   

We also define
\begin{align}\label{eq:theta}
  \theta = \theta(\hat{x}, \hat{x}_0, \alpha) := \frac{\hat{x}-\alpha \hat{x}_0}{|\hat{x}-\alpha \hat{x}_0|}. 
\end{align}

\begin{theorem}\label{thm:main}
Let $\psi_0$ and $\psi_1$ be as in Problem \ref{prbl:1}, where $A = A(|x_0|),$ and $\kappa$ is fixed. We consider $x,$ $x_0,$ $y \in \mathbb{R}^d$  as in \eqref{eq:notations}, where 
\begin{align}
&\sin\left(\kappa(\theta-\hat{x}, \zeta)\right) \neq 0, \label{eq:2.10.1} \\
    &\text{either } \hat{x} \neq \hat{x}_0, \quad \text{or} \quad \alpha\neq 1. \label{eq:2.11.1}
\end{align}
Let $f_1$ be the far-field pattern for $\psi_1,$ defined in \eqref{eq:psif1}.
Then $f_1(\hat{x})$ is approximately determined by the intensities $|\psi_0(x)+\psi_1(x)|^2$ and $|\psi_0(y)+\psi_1(y)|^2$  via the following formula:
\begin{align}
\begin{aligned}\label{eq:2.6}
    &f_1(\hat{x}) = f_{1,1}(x, x_0, y) + \frac{1}{D}\left({\cal O}\left(r^{-1}\right)+{\cal O}\left(\frac{1}{A(\alpha r)}\right) \right), \text{ as } r  \to +\infty, \quad \zeta \text{ is fixed,} \\
    &f_{1,1}(x, x_0, y) := \frac{e^{i\varphi}}{D} \left(e^{i\kappa|y-x_0|-i\kappa|y|}b(x, x_0) -e^{i\kappa|x-x_0|-i\kappa|x|}b(y, x_0)\right), \\
    &D = 2i \sin  \kappa(|y-x_0|-|x-x_0|+|x|-|y| ),     
\end{aligned}
\end{align}
where $r$ is the parameter in  \eqref{eq:notations},  $b$ is defined by \eqref{eq:b}, and $\varphi$ arises in \eqref{eq:2.5}.

\end{theorem}

We consider \eqref{eq:2.6} assuming that $D\neq 0$. In addition, 
\begin{align}\label{eq:2.13}
    D = 2i \sin \left[\kappa(\theta-\hat{x}, \zeta) +{\cal O}\left(r^{-1}\right)\right], \quad r\to +\infty, \quad \zeta \text{ is fixed}.    
    \end{align}

Note also that we consider  \eqref{eq:2.6} assuming that
$A(\alpha r)$ is sufficiently large, as $r\to+\infty$.

Theorem \ref{thm:main} and formula \eqref{eq:2.13} are proved in Section \ref{sec:3}, where more detailed versions of ${\cal O}$ in \eqref{eq:2.6}, \eqref{eq:2.13} are also given.

Next, formulas \eqref{eq:2.6}, \eqref{eq:2.13} admit direct applications to Problem \ref{prbl:1} for the case when ${\cal U}' = {\cal U}'' = X,$ where $X$ is a hyperplane in ${\cal U},$ and the distance between $X$ and the origin $\{0\}\in \mathbb{R}^d$ is sufficiently large.

Let 
\begin{align}
&X = X_{s, \omega} = \{ x \in \mathbb{R}^d: \quad (x, \omega) = s\}, \quad s\geq 0, \quad \omega \in \mathbb{S}^{d-1}, \label{eq:X}\\
&{\mathbb S}^+_{\omega} = \{\hat{x} \in \mathbb{S}^{d-1}: (\hat{x}, \omega)>0\}, \quad \omega \in \mathbb{S}^{d-1}. \label{eq:S}
\end{align}

Using Theorem \ref{thm:main} and formula \eqref{eq:psif1} we obtain the following result.
\begin{corollary}\label{cor:1}
Let $\psi_0,$ $\psi_1$ be as in Theorem \ref{thm:main}.  Let $\omega, \hat{x}_0 \in \mathbb{S}^{d-1},$ $\alpha>0.$ 
Then the following  formula holds for approximate finding  $\psi_1$ on $X_{s, \omega}$ from $|\psi_0+\psi_1|^2$ on $X_{s, \omega}$, where $s$ is sufficiently large:
\begin{align}
\begin{aligned}\label{eq:2.17}
&\psi_1(x)  = \frac{e^{i\kappa r}}{r^{\frac{d-1}{2}}}f_{1, 1}\left(x, x_0, y\right)+ \frac{1}{r^{\frac{d-1}{2}}D}\left({\cal O}\left(r^{-1}\right)+{\cal O}\left(\frac{1}{A(\alpha r)}\right) \right) +{\cal O}\left(\frac{1}{r^{\frac{d+1}{2}}}\right),  \quad \text{ as } r\to +\infty, \\
&x \in X_{s, \omega},  \quad y = x+\zeta \in X_{s, \omega}, \quad |x| = r, \quad x = r \hat{x}, \quad |x_0| = \alpha r, \quad x_0 = \alpha r \hat{x}_0,  
\end{aligned}
\end{align}
      where $f_{1, 1}$ and $D$ are defined as in \eqref{eq:2.6}, and formula \eqref{eq:2.17} is considered for each fixed pair $\hat{x} \in \mathbb{S}^+_{\omega}$ and $\zeta \in X_{0, \omega}$ under  conditions \eqref{eq:2.10.1}, \eqref{eq:2.11.1}.  
      
\end{corollary}

We consider \eqref{eq:2.17} assuming that $D\neq 0.$ In addition, in view of \eqref{eq:2.13}, there is a large choice of $\zeta\in X_{0, \omega}$  such that
\begin{align}\label{eq:2.18}
    \lim_{r\to +\infty} D = 2i \sin \kappa(\theta-\hat{x}, \zeta) \neq 0,
\end{align}
\noindent under the condition that  
\begin{align}\label{eq:2.19}
\pi_{\omega}(\theta-\hat{x})\neq 0,
\end{align}
where
\begin{align}
    &\pi_{\omega} \text{ is the orthogonal projector of $\mathbb{R}^d$ on } X_{0, \omega},
\end{align}
$\theta$ is defined in \eqref{eq:theta}. 

In addition, first of all, it is natural to consider  Corollary \ref{cor:1} assuming that $\hat{x}_0 \in -\mathbb{S}^+_{\omega},$ where $\mathbb{S}^+_{\omega}$ is defined by \eqref{eq:S}. In this case the following lemma holds.
\begin{lemma}\label{lem:2} Let $\omega \in \mathbb{S}^{d-1},$ $\hat{x}_0 \in -\mathbb{S}^+_{\omega},$ $\alpha>0.$  Then condition \eqref{eq:2.19} holds if 
\begin{align}\label{eq:xse}
\hat{x} \in \mathbb{S}^+_{\omega}\setminus \hat{{\cal E}}_{\omega, \hat{x}_0, \alpha},   
\end{align}
where $\hat{{\cal E}}_{\omega, \hat{x}_0, \alpha}$ is an exceptional subset of $\mathbb{S}^+_{\omega},$   such that 
\begin{align}
   \hat{{\cal E}}_{\omega, \hat{x}_0, \alpha} \subset span(\{\omega, \hat{x}_0\}), \quad \#\hat{{\cal E}}_{\omega, \hat{x}_0, \alpha} \leq 6, 
\end{align}
 where $\#$ denotes the number of elements.

\end{lemma}

Formula \eqref{eq:2.17} follows just from Theorem \ref{thm:main} and formula \eqref{eq:psif1}. Lemma \ref{lem:2} is proved in Section \ref{sec:4}.

\begin{remark}\label{rem:3}
In fact, formula \eqref{eq:2.17} is mainly of interest for $\hat{x}\in \mathbb{S}^+_{\omega}$  and $\zeta \in X_{0, \omega},$ when the limit in \eqref{eq:2.18} is not zero. The point is that, 
for the case when $\hat{x}_0\in -\mathbb{S}^+_{\omega},$ this limit is non-zero for any $\hat{x}$ satisfying \eqref{eq:xse} with a large choice of appropriate $\zeta \in X_{0, \omega}$. Therefore,  
formula \eqref{eq:2.17} can be used efficiently for  $x \in X_{s, \omega}\setminus {\cal E}_{\omega, \hat{x}_0, \alpha},$ where $s$ is sufficiently large, and
\begin{align}
{\cal E}_{\omega, \hat{x}_0, \alpha} := \{x \in X_{s, \omega}: \hat{x} \in \hat{\cal{E}}_{\omega, \hat{x}_0, \alpha} \}; 
\end{align}
in addition, $\#{\cal E}_{\omega, \hat{x}_0, \alpha} = \#\hat{\cal{E}}_{\omega, \hat{x}_0, \alpha}.$  
In formula  \eqref{eq:2.17} we also assume that
$A(\alpha r)$ is sufficiently large, as $r\to+\infty$.
\end{remark}

 For the case when $\psi_0$ is a plane-wave solution  of the Helmholtz equation \eqref{eq:schrod},  prototypes of Theorem \ref{thm:main} are given in \cite{N2015}, \cite{NS2026}, prototypes of Corollary \ref{cor:1} are given in  \cite{NS2026}.


\section{Proof of Theorem \ref{thm:main} and formula \eqref{eq:2.13}}\label{sec:3}

\noindent Let 
\begin{align}\label{eq:psi01}
    \psi = \psi_0+\psi_1.
\end{align}

\noindent Using  \eqref{eq:psi_green}, \eqref{eq:psif1}, \eqref{eq:psi01}, we get
\begin{align}
\begin{aligned}
 &|\psi|^2 = \psi \overline{\psi} =  |\psi_0(x)|^2 +\frac{\overline{A'}e^{-i\kappa|x-x_0|} \psi_1(x)}{|x-x_0|^{(d-1)/2}}+\frac{A'e^{i\kappa|x-x_0|} \overline{\psi_1(x)}}{|x-x_0|^{(d-1)/2}}  +\psi_1(x) \overline{\psi_1(x)} +{\cal O}\left(\frac{|A'|}{|x|^{(d-1)/2}|x-x_0|^{(d+1)/2}}\right) = \\
     & =|\psi_0(x)|^2 +\frac{\overline{A'}e^{-i\kappa|x-x_0|} e^{i\kappa |x|} f_1(\hat{x})}{|x|^{(d-1)/2}|x-x_0|^{(d-1)/2}}+ \frac{A'e^{i\kappa|x-x_0|-i\kappa|x|} \overline{f_1(\hat{x})}}{|x|^{(d-1)/2}|x-x_0|^{(d-1)/2}}  + \\
&+\frac{|f_1(\hat{x})|^2}{|x|^{d-1}} +  {\cal O}\left(\frac{1}{|x|^{d}}\right)+{\cal O}\left(\frac{|A'|}{|x|^{(d+1)/2}|x-x_0|^{(d-1)/2}}\right) +{\cal O}\left(\frac{|A'|}{|x|^{(d-1)/2}|x-x_0|^{(d+1)/2}}\right),     
\end{aligned}\label{eq:3.4}
\end{align}
as $|x| \to +\infty$ and $|x-x_0| \to +\infty.$


In view of \eqref{eq:b}, \eqref{eq:3.4}, we have that 
\begin{align}
&\begin{aligned}
b(x, x_0) = &e^{-i(\varphi+\kappa|x-x_0|-\kappa|x|)} f_1(\hat{x}) + e^{i(\varphi+\kappa|x-x_0|-\kappa |x|)} \overline{f_1(\hat{x})}+ \frac{|f_1(\hat{x})|^2 |x-x_0|^{(d-1)/2}}{|A'||x|^{(d-1)/2}} + \\
&+{\cal O}\left(\frac{|x-x_0|^{(d-1)/2}}{|A'||x|^{(d+1)/2}}\right)+{\cal O}\left(\frac{1}{|x|}\right) +{\cal O}\left(\frac{1}{|x-x_0|}\right), \quad \text{ as $|x|\to +\infty$ and $|x-x_0|\to +\infty;$ }  
\end{aligned}\label{eq:3.5} \\
&\begin{aligned}
b(y, x_0) = &e^{-i(\varphi+\kappa|y-x_0|-\kappa|y|)} f_1(\hat{y}) + e^{i(\varphi+\kappa|y-x_0|-\kappa |y|)} \overline{f_1(\hat{y})}+ \frac{|f_1(\hat{y})|^2 |y-x_0|^{(d-1)/2}}{|A'||y|^{(d-1)/2}} + \\
&+{\cal O}\left(\frac{|y-x_0|^{(d-1)/2}}{|A'||y|^{(d+1)/2}}\right)+{\cal O}\left(\frac{1}{|y|}\right) +{\cal O}\left(\frac{1}{|y-x_0|}\right), \quad \text{ as $|y|\to +\infty$ and $|y-x_0|\to +\infty.$ }     
\end{aligned}\label{eq:3.6}
\end{align}

Using \eqref{eq:notations}, \eqref{eq:3.5}, and \eqref{eq:3.6}, we get
\begin{align}
&\begin{aligned}
b(x, x_0) = &e^{-i(\varphi+\kappa|x-x_0|-\kappa|x|)} f_1(\hat{x}) + e^{i(\varphi+\kappa|x-x_0|-\kappa |x|)} \overline{f_1(\hat{x})}+ {\cal O}\left(r^{-1}\right)+{\cal O}\left(\frac{1}{A(\alpha r)}\right) ,  
\end{aligned}\label{eq:3.7} \\
&\begin{aligned}
b(y, x_0) = &e^{-i(\varphi+\kappa|y-x_0|-\kappa|y|)} f_1(\hat{x}) + e^{i(\varphi+\kappa|y-x_0|-\kappa |y|)} \overline{f_1(\hat{x})}+ {\cal O}\left(r^{-1}\right)+{\cal O}\left(\frac{1}{A(\alpha r)}\right),    
\end{aligned}\label{eq:3.8}
\end{align}
as $r\to +\infty.$

Note that in \eqref{eq:3.8} we also used that 
\begin{align}
    f_1(\hat{y}) = f_1(\hat{x}) + {\cal O}\left(\frac{|\zeta|}{r}\right),  \text{ as } r\to+\infty,
\end{align}
since  $f_1$ is smooth on $\mathbb{S}^{d-1}$ (see, for example, \cite{W}, \cite{Karp61}), and
\begin{align}
    \hat{y} = \hat{x} + {\cal O}\left(\frac{|\zeta|}{r}\right), \text{ as } r\to+\infty. 
\end{align}

Theorem \ref{thm:main} follows from the approximate 2 × 2  linear system \eqref{eq:3.7}, \eqref{eq:3.8} for $f_1$, $\overline{f_1}$. In particular, $D$ in \eqref{eq:2.6} is the
determinant of the corresponding  matrix, arising in this system; and it is assumed that $D \neq 0.$

Next, to obtain formula \eqref{eq:2.13}, we use the formula for $D$ in \eqref{eq:2.6} and the following asymptotic formulas
\begin{align}
&\begin{aligned}
&|y| = |x+\zeta| = |x|\left(\hat{x} +\frac{\zeta}{|x|}, \hat{x}+\frac{\zeta}{|x|}\right)^{1/2} = |x|\left( 1 +2\frac{(\hat{x}, \zeta)}{|x|}+{\cal O}\left(\frac{|\zeta|^2}{|x|^2}\right)\right)^{1/2} = \\
&=|x|\left(1+ \frac{(\hat{x}, \zeta)}{|x|}+{\cal O}\left(\frac{|\zeta|^2}{|x|^2}\right)\right), \quad |x|\to +\infty,    
\end{aligned} \\
&|y-x_0| = |x-x_0+\zeta| = |x-x_0| \left(1 + \frac{(\widehat{x-x_0}, \zeta)}{|x-x_0|}+{\cal O}\left(\frac{|\zeta|^2}{|x-x_0|^2}\right)\right), \quad |x-x_0|\to +\infty.
\end{align}
Using \eqref{eq:notations} and \eqref{eq:theta}, we obtain
\begin{align}
\begin{aligned}\label{eq:3.11}
&|y-x_0|-|x-x_0|+|x|-|y| = \\
&=|x-x_0|\left(1+ \frac{(\widehat{x-x_0}, \zeta)}{|x-x_0|}+{\cal O}\left(\frac{|\zeta|^2}{|x-x_0|^2}\right)\right) -  |x-x_0| + |x| - |x|\left(1+ \frac{(\widehat{x}, \zeta)}{|x|}+{\cal O}\left(\frac{|\zeta|^2}{|x|^2}\right)\right) = \\
&=(\theta-\hat{x}, \zeta) +{\cal O}\left(|\zeta|^2(\frac{1}{|x|}+\frac{1}{|x-x_0|})\right), \quad \text{as } |x|\to+\infty \text{ and } |x-x_0|\to+\infty.     
\end{aligned} 
\end{align}
In addition, in view of \eqref{eq:notations}, 
\begin{align}\label{eq:3.12}
    \frac{1}{|x|} +\frac{1}{|x-x_0|} = \frac{1}{r} \left(1+\frac{1}{|\hat{x}-\alpha \hat{x}_0|}\right). 
\end{align}
Finally, formula \eqref{eq:2.13} follows from formula for $D$ in \eqref{eq:2.6}  and formulas \eqref{eq:3.11},  \eqref{eq:3.12}, \eqref{eq:2.11.1}.

\section{Proof of Lemma \ref{lem:2}}\label{sec:4}




\noindent Let
\begin{align}
\begin{aligned}\label{eq:2.16}
&\hat{x}_0 = -\omega \xi_1+ \gamma \xi_2, \quad \xi = (\xi_1, \xi_2) \in \mathbb{S}^1, \quad \xi_1>0,\\
&\hat{x} = \omega  \eta_1+\gamma \eta_2 , \text{ for } d=2, \quad \eta = (\eta_1, \eta_2) \in \mathbb{S}^1 , \quad\eta_1>0, \\ &\hat{x} = \omega  \eta_1+\gamma \eta_2 +\nu \eta_3, \text{ for } d\geq 3, \quad \eta = (\eta_1, \eta_2, \eta_3) \in \mathbb{S}^2, \quad \eta_1>0, \\
&\omega, \gamma, \nu \in \mathbb{ S}^{d-1}, \quad (\gamma, \omega) = 0, \quad  (\nu, \omega)= (\nu, \gamma) = 0,
\end{aligned}
\end{align}
where the orthonormal vectors $\omega, \gamma, \nu$  are fixed.

In view of \eqref{eq:theta}, \eqref{eq:2.16},  we have that:
\begin{align}
&\theta-\hat{x} = \frac{\omega(\eta_1+\alpha \xi_1)+\gamma(\eta_2-\alpha \xi_2)}{\sqrt{(\eta_1+\alpha \xi_1)^2+(\eta_2-\alpha \xi_2)^2}}- \omega  \eta_1-\gamma \eta_2, \label{eq:4.2} \\
&\pi_{\omega}(\theta-\hat{x}) = \frac{\gamma(\eta_2-\alpha \xi_2)}{\sqrt{(\eta_1+\alpha \xi_1)^2+(\eta_2-\alpha \xi_2)^2}}-\gamma \eta_2, \label{eq:4.3}
\end{align}
where $d=2;$ and 
\begin{align}
&\theta-\hat{x} = \frac{\omega(\eta_1+\alpha \xi_1)+\gamma(\eta_2-\alpha \xi_2)+\nu \eta_3}{\sqrt{(\eta_1+\alpha \xi_1)^2+(\eta_2-\alpha \xi_2)^2+\eta_3^2}}- \omega  \eta_1-\gamma \eta_2 -\nu \eta_3, \\
&\begin{aligned}\label{eq:pi_omega_form}
&\pi_{\omega}(\theta-\hat{x}) = \frac{\gamma(\eta_2-\alpha \xi_2)+\nu \eta_3}{\sqrt{(\eta_1+\alpha \xi_1)^2+(\eta_2-\alpha \xi_2)^2+\eta_3^2}}-\gamma \eta_2 -\nu \eta_3 = \\
&=\gamma\left(\frac{(\eta_2-\alpha \xi_2)}{\sqrt{(\eta_1+\alpha \xi_1)^2+(\eta_2-\alpha \xi_2)^2+\eta_3^2}}- \eta_2 \right) +\nu \eta_3 \left(\frac{1}{\sqrt{(\eta_1+\alpha \xi_1)^2+(\eta_2-\alpha \xi_2)^2+\eta_3^2}}-1 \right),    
\end{aligned}
\end{align}
where $d\geq 3.$

Next, for $d\geq 3, $ we prove that 
\begin{align}\label{eq:4.6}
    \text{if $\eta_3\neq 0$ in \eqref{eq:2.16}, then condition \eqref{eq:2.19} holds.}
\end{align}

Using \eqref{eq:pi_omega_form}, where $\eta_3\neq 0,$ we have that 
\begin{align}
    \pi_{\omega}(\theta-\hat{x})= 0 \quad \text{ iff} \quad \begin{cases}\label{eq:4.7}
&\eta_3 \left(\frac{1}{\sqrt{(\eta_1+\alpha \xi_1)^2+(\eta_2-\alpha \xi_2)^2+\eta_3^2}}-1 \right) = 0, \\  
&\frac{(\eta_2-\alpha \xi_2)}{\sqrt{(\eta_1+\alpha \xi_1)^2+(\eta_2-\alpha \xi_2)^2+\eta_3^2}}- \eta_2  =0.    
\end{cases}
\end{align}


Using  the presentations for $\hat{x}_0,$ $\hat{x}$  in \eqref{eq:2.16}, and using the system of equations for $\eta$ in  \eqref{eq:4.7}, we get 
\begin{align}
\begin{cases}
&\xi_1^2+\xi_2^2 = 1, \quad \xi_1>0,\\
&\eta_1^2+\eta_2^2+\eta_3^2 = 1, \quad \eta_1>0,  \\
&(\eta_1+\alpha \xi_1)^2+(\eta_2-\alpha \xi_2)^2+\eta_3^2 = 1, \\
&\alpha \xi_2 = 0,
\end{cases}
\end{align}
that is,
\begin{align}
\begin{cases}\label{eq:4.9}
&\eta_1^2+\eta_2^2+\eta_3^2 = 1, \quad \eta_1>0, \\
&(\eta_1+\alpha \xi_1)^2+\eta_2^2+\eta_3^2 = 1, \\
&\xi_1 =  1.
\end{cases}
\end{align}
However, using also that $\alpha>0,$ one can see that system \eqref{eq:4.9} has no solution with respect to $\eta\in \mathbb{S}^2.$

Thus, \eqref{eq:4.6} is proved.

Next, for $d=2$ and for $d\geq 3,$ where $\eta_3=0,$ we prove that 
\begin{align}\label{eq:4.10}
    \text{if $\hat{x}_0=-\omega,$ then condition \eqref{eq:2.19} holds if $\hat{x} \neq \omega.$}
\end{align}
In this special case, in view of \eqref{eq:4.2}--\eqref{eq:pi_omega_form}, we have that 
\begin{align}\label{eq:4.11}
    \pi_{\omega}(\theta-\hat{x}) = \gamma\left(\frac{\eta_2}{\sqrt{(\eta_1+\alpha )^2+\eta^2_2}}- \eta_2 \right).  
\end{align}
Using that $\eta = (\eta_1, \eta_2) \in \mathbb{S}^1,$ $\eta_1>0,$ and $\alpha>0,$ one can see that 
\begin{align}\label{eq:4.12}
\sqrt{(\eta_1+\alpha )^2+\eta^2_2}>1.
\end{align}
In view of \eqref{eq:4.11}, \eqref{eq:4.12}, one can see  that  $\pi_{\omega}(\theta-\hat{x})\neq 0$ if $\eta_2 \neq 0,$ i.e., if $\hat{x}\neq \omega.$ 

Thus, \eqref{eq:4.10} is proved.

Formulas \eqref{eq:2.16} and properties \eqref{eq:4.6}, \eqref{eq:4.10} imply that 
\begin{align}
    \hat{{\cal E}}_{\omega, \hat{x}_0, \alpha} \subset span(\{\omega, \hat{x}_0\}).
\end{align}

Finally, we prove that \begin{align}\label{eq:ne}
    \#\hat{{\cal E}}_{\omega, \hat{x}_0, \alpha} \leq 6.
\end{align}

In view of \eqref{eq:4.6}, it is sufficient to consider the case $d=2,$ or the case $d\geq 3,$ where $\eta_3 = 0.$ In this case, 
in view of \eqref{eq:4.3}, \eqref{eq:pi_omega_form}, we have that \begin{align}\label{eq:iff}
    \pi_{\omega}(\theta-\hat{x}) = 0 \quad \text{ iff } \quad \frac{(\eta_2-\alpha \xi_2)^2}{(\eta_1+\alpha \xi_1)^2+(\eta_2-\alpha \xi_2)^2} = \eta_2^2. 
\end{align}
In turn, \eqref{eq:iff} can be rewritten as,
\begin{align}
\begin{aligned}
&(\eta_2-\alpha \xi_2)^2 = \eta_2^2-2\alpha \eta_2\xi_2+\alpha^2\xi_2^2= \eta_2^2((\eta_1+\alpha \xi_1)^2+(\eta_2-\alpha \xi_2)^2)= \\
&=\eta_2^2\left(1+2\alpha \eta_1\xi_1+\alpha^2\xi_1^2-2\alpha \eta_2 \xi_2 +\alpha^2\xi_2^2\right) = \eta_2^2\left(1+\alpha^2+2\alpha \eta_1\xi_1-2\alpha \eta_2 \xi_2\right)= \\
&=\eta_2^2+\alpha \eta_2^2\left(\alpha+2 \eta_1\xi_1-2 \eta_2 \xi_2\right),     
\end{aligned}
\end{align}
and further as
\begin{align}\label{eq:4.17}
&\alpha \xi_2^2 = \eta_2^2\left(\alpha+2 \eta_1\xi_1\right)+2 \eta_2 \xi_2(1-\eta_2^2) = \eta_2^2\left(\alpha+2 \eta_1\xi_1\right)+2 \eta_1^2 \eta_2 \xi_2.
\end{align}
To study \eqref{eq:iff} with respect to $\eta\in \mathbb{S}^1,$ we rewrite \eqref{eq:4.17} in terms of $\lambda = \eta_1+i\eta_2,$ where 
\begin{align}
&\eta_1 = (\lambda+\lambda^{-1})/2, \quad \eta_2 = (\lambda-\lambda^{-1})/(2i). 
\end{align}
We get
\begin{align}
&\begin{aligned}
&\alpha \xi_2^2 = (\lambda^2-2+\lambda^{-2})\frac{-1}{4} (\alpha+\lambda \xi_1+\lambda^{-1}\xi_1)+ \frac{1}{2}(\lambda^2+2+\lambda^{-2})\frac{\lambda-\lambda^{-1}}{2i} \xi_2 = \\
&=\frac{1}{4\lambda^3}\left[-(\lambda^4-2\lambda^2+1) (\alpha \lambda+\lambda^2 \xi_1+\xi_1)-i (\lambda^4+2\lambda^2+1)(\lambda^2-1)\xi_2\right] ,    
\end{aligned}
\end{align}
and finally
\begin{align}\label{eq:4.20}
&\lambda^6(-\xi_1-i\xi_2)+\lambda^5(-\alpha)+\lambda^4(\xi_1-i\xi_2)+\lambda^3(-4\alpha \xi_2^2+2\alpha)+\lambda^2(\xi_1+i\xi_2)+\lambda (-\alpha) + (-\xi_1+i\xi_2)=0,
\end{align}
where $\xi = (\xi_1, \xi_2)\in \mathbb{S}^1.$

Formula \eqref{eq:ne} follows from \eqref{eq:iff} rewritten in the form \eqref{eq:4.20}. This completes the proof of Lemma \ref{lem:2}.

Note also that if $\lambda \in \mathbb{C}$ satisfy \eqref{eq:4.20}, then $\overline{\lambda^{-1}}$ also satisfies \eqref{eq:4.20}.

\section{Acknowledgements}

The research was supported by Russian Science Foundation, grant № 25-71-00116, https://rscf.ru/project/25-71-00116/ (V.N. Sivkin).

\bibliographystyle{alpha}

\vskip 3mm

Roman G. Novikov

CMAP, CNRS, \'Ecole polytechnique, Institut Polytechnique de Paris, 91128 Palaiseau, France 


E-mail: roman.novikov@polytechnique.edu
\vskip 3mm

Vladimir N. Sivkin 

HSE University, Moscow, Russia

E-mail: sivkin96@yandex.ru

\end{document}